\documentclass{amsart}

\usepackage{amssymb,url}

\usepackage{ytableau}

\newtheorem{theorem}{Theorem}[section]

\theoremstyle{definition}
\newtheorem{definition}[theorem]{Definition}
\newtheorem{example}[theorem]{Example}

\theoremstyle{remark}
\newtheorem{remark}[theorem]{Remark}

\numberwithin{equation}{section}
\newcommand{\complexi}{\boldsymbol{i}}
\newcommand{\lb}{\left\{ }
\newcommand{\rb}{\right\} }

\newcommand{\lp}{\left(}
\newcommand{\rp}{\right)}

\newcommand{\series}{\sum^\infty}
\newcommand{\complex}{\mathbb{C}}

\newcommand{\nats}{\mathbb{N}}

\newcommand{\ints}{\mathbb{Z}}

\newcommand{\Ii}{\mathcal{I}}

\newcommand{\Kk}{\mathcal{K}}

\newcommand{\BPsi}{\boldsymbol{\Psi}}
\newtheorem{innercustomgeneric}{\customgenericname}
\providecommand{\customgenericname}{}
\newcommand{\newcustomtheorem}[2]{%
  \newenvironment{#1}[1]
  {%
   \renewcommand\customgenericname{#2}%
   \renewcommand\theinnercustomgeneric{##1}%
   \innercustomgeneric
  }
  {\endinnercustomgeneric}
}

\newcustomtheorem{customthm}{Theorem}
\newcustomtheorem{customlemma}{Lemma}

\begin{document}


\title[Exceptional Hermites and CM Pairs]{Exceptional Hermite Polynomials and\\ Calogero-Moser Pairs
}


\author{Luke Paluso}
\address{Department of Mathematics / College of Charleston /
  Charleston SC 29401}
\email{palusol@g.cofc.edu}

\author{Alex Kasman}
\address{Department of Mathematics / College of Charleston /
  Charleston SC 29401}
\email{kasmana@cofc.edu}


\begin{abstract}
There are two equivalent descriptions of George Wilson's adelic
Grassmannian $Gr^{ad}$, one in terms of differential ``conditions''
and another in terms of Calogero-Moser Pairs.  The former approach was used in the
2020 paper by Kasman-Milson which found that each family of
Exceptional Hermite Polynomials has a generating function which lives
in $Gr^{ad}$.  This suggests that Calogero-Moser Pairs should also
be useful in the study of Exceptional Hermite Polynomials, but no
researchers have pursued that line of inquiry prior
to the first author's thesis.  
The purpose of this note is to summarize highlights from that
thesis, including a novel formula for Exceptional Hermite Polynomials
in terms of Calogero-Moser Pairs and a theorem utilizing this
correspondence to produce explicit finitely-supported distributions
which annihilate them.

\bigskip

{\scriptsize\noindent
(First published as \textit{Contemporary Mathematics}
Vol.~822 pp.~167--180\\
https://doi.org/10.1090/conm/822/16481.  \copyright~2025 Amer.~Math.~Society.)}
\end{abstract}

\maketitle
\section{Introduction}
The mathematical objects mentioned in the title have separately been the subjects of many
papers.  Section~\ref{sec2} below will provide additional details, but for the sake of this introduction let us just briefly recall:
\begin{itemize}
\item  A pair of $N\times N$ constant complex matrices $(X,Z)$
  whose commutator differs from the identity by a rank one matrix are
  called a Calogero-Moser Pair (CM Pair).  Interest in them first arose in the
  1970s to explain the integrable dynamics of the Calogero-Moser particle
  system using Hamiltonian reduction \cite{KKSCMPAIRS}. In 1998, Wilson showed that the rank one bispectral wave functions
  of the KP Hierarchy can all be written in terms of CM Pairs.
 \cite{Wilson1998}.
  \item
Exceptional orthogonal polynomials, a
generalization of the classical orthogonal polynomials sharing many of
their important properties despite lacking polynomials of certain
degrees, were first defined in 2009 \cite{UllateKamranMilson2009}
and the particular example of the
Exceptional Hermite Polynomials (XHPs) was first considered in 2014
\cite{UllateGrandatiMilsonXHPS}.   
\end{itemize}
It is not clear from those descriptions that there is any meaningful
connection between XHPs and CM Pairs and they have not appeared
together in any prior papers.  However, it was shown in
2020 that each family of XHPs has a generating function which is a
rank one bispectral wave function of the KP Hierarchy
\cite{MilsonKasman2020}.  Consequently, it follows from Wilson's
observation that the generating functions for XHPs can be written in
terms of CM Pairs.

Finding a formula for XHPs in terms of CM Pairs and applications
for the correspondence between them
was the thesis project assigned to the first
author by the second author.  This note constitutes a research
announcement aiming to convey the main results found in
\cite{Thesis}.  The thesis and its references should be consulted for
details and complete proofs.

\section{Background Material}\label{sec2}

\subsection{Partitions and Young Diagrams}\label{subsection:Partitions and Young Diagrams}\label{sec:frob}
For a fixed natural number
$N \in \nats$, a \textit{partition of $N$} is a non-increasing sequence of non-negative integers $\lambda=\left\{ \lambda_n \right\}_{n \in \nats}$
such that 
\begin{align*}
    |\lambda|= \series_{n=1} \lambda_n =N.
\end{align*}
Since the elements of the sequence must eventually be zero, a
partition can be equivalently written as a finite list by omitting all
of the $0$ entries.  The length of the partition $\ell=\ell(\lambda)$,
is the number of non-zero entries in the sequence.

For example, the partition $\lambda=\{4,3,2,0,0,\ldots\}$ is a
partition of $|\lambda|=9$ which has length $\ell(\lambda)=3$ and
hence could also be written $\lambda=\{4,3,2\}$.

Following \cite{MilsonKasman2020} we define some sets of integers
corresponding to the choice of a partition $\lambda$,
\begin{align}
    \mathcal{M}^{(\lambda)}&=\left\{ m_i = \lambda_i-i :i\in \nats \right\} \label{eqn:MayaDiagramDef}\\
   \mathcal{K}_n^{(\lambda)}&=\left\{ m_i+n : 1\leq i \leq n \right\} \label{eqn:indexSet} \\ 
    \mathcal{I}^{(\lambda)}&=\ints \setminus \mathcal{M}^{(\lambda)}+|\lambda| = \lb z+|\lambda|\; : \; z\in \ints \setminus \mathcal{M}^{(\lambda)}  \rb \label{eqn:scriptI}.
\end{align}
$\mathcal{M}^{(\lambda)}$ is known as the \emph{Maya Diagram} of a
partition $\lambda$. The set defined in \eqref{eqn:indexSet} depends on the
choice of an integer $n$ and will be of special interest in the
case that $n=\ell(\lambda)$ is the length of the partition: $\mathcal{K_\ell}^{(\lambda)}$ is known as the
\emph{Index Set} of a partition and we will just write it as
$\Kk^{(\lambda)}=\lb k_1,\dots,k_\ell \rb$.  It is also noteworthy that $\Ii^{(\lambda)} \subset \nats$.

A \emph{Young Diagram} is a finite array of boxes arranged in
left-justified rows, with the rows in non-increasing order according
to their lengths. There is a natural bijection between partitions and
Young Diagrams which associates to
$\lambda=\{\lambda_1,\lambda_2,\ldots\}$ the diagram whose $i^{th}$
row has $\lambda_i$ boxes. 

An \textit{alternate} method for representing partitions is the
\textit{Frobenius notation}.  Consider the Young Diagram corresponding
to a partition $\lambda$.  The number of
boxes to the right of the box in the $i^{th}$ row and $i^{th}$ column in
the Young Diagram is denoted $\alpha_i$. Similarly, $\beta_i$ is the
number of boxes below the box in the $i^{th}$ row and $i^{th}$
column. The partition corresponding to that Young Diagram can then be denoted
$\lambda=\left\{ (\alpha_1,\dots, \alpha_{N'})|(\beta_1,\dots,
  \beta_{N'}) \right\}$ where $N'$ is the highest integer $j$ so that $\lambda_j \geq j$ and where the Young Diagram has a box in position $(j,j)$.

 The Young Diagram for the partition $\lambda={ 4,3,2 }$ is shown in Figure \ref{fig:4,3,2}.
\begin{figure}
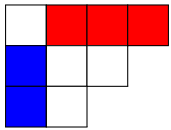

\begin{align*}
    \begin{ytableau}
        *(white) & *(red) & *(red) &*(red) \\
        *(blue) & *(white) & *(white) \\
        *(blue) & *(white)
    \end{ytableau}
\end{align*}
    \caption{$\lambda=\{4,3,2\}=\lb (3,1) | (2,1) \rb$}
    \label{fig:4,3,2}
\end{figure}
All boxes to the right of the box in the first column and first row
are colored red corresponding with $\alpha_1=3$. Similarly the number
of boxes below the box in the first column and first row are colored
blue, corresponding with $\beta_1=2$. Similarly we can conclude
$\alpha_2=\beta_2=1$. Hence, in Frobenius notation the partition
corresponding to that Young Diagram is $\lambda=\left\{ (3 , 1) \; | \; (2 , 1)\right\}$.

\subsection{Bivariate Exceptional Hermite Polynomials}\label{subsection:XHPS}
Before we get to the XHPs, we
begin by defining the classical Hermite polynomials.  Readers may
already be familiar with them and inclined to skip over this definition, but note that here we work with them in an
unusual \textit{bivariate} form.

For the purposes of this paper, the Hermite polynomials $H_n(x,y)$ are defined by the
generating function 
$$e^{xz+yz^2}=\series_{n=0}H_n\frac{z^n}{n!}.$$
In the special case $y=-1/4$ they yield the standard univariate
Hermite polynomials.  Moreover, for any negative value of $y$ they are
orthogonal satisfying
 \begin{align*}
     \langle H_{n},H_{m} \rangle_{H_y} = 2(-\pi y)^{\frac{1}{2}}(-2y)^n n!\delta_{n,m}.
 \end{align*}
with respect to the inner product
\begin{align*}
    \langle f(x),g(x) \rangle_{H_y} = \int_{-\infty}^\infty f(x)g(x) e^{\frac{x^2}{4y}}dx.
\end{align*}

To each partition $\lambda$, we
associate a family of bivariate Exceptional Hermite Polynomials
(briefly referred to as the $\lambda$-XHPs) \cite{MilsonKasman2020}.  This is the set
$\left\{ H_{n}^{(\lambda)}\ :\ n\in \Ii^{(\lambda)} \right\}$
of polynomials defined by the formula
\begin{align}
  & H_n^{(\lambda)}= \frac{\text{Wr}(H_{k_\ell},H_{k_{\ell-1}}\dots,H_{k_1},H_{n-N+\ell})}{\prod_{i<j}(k_i-k_j)\prod_{i}(n-N+\ell-k_i)}& &\text{for}\;\;\; n \in \Ii^{(\lambda)},& \label{eqn:defofXHPS}
\end{align}
where the Wronskian involves differentiation with respect to the
variable $x$, $N=|\lambda|$, 
and where
$\Ii^{(\lambda)}$ and index set
$\Kk^{(\lambda)}=\{k_1,\ldots,k_\ell\}$ are as in Section
\ref{subsection:Partitions and Young Diagrams}.

The polynomial $H_n^{(\lambda)}(x,y)$ is homogeneous of weighted
degree $n$ where $\text{deg}(x)=1$ and
$\text{deg}(y)=2$. Consequently, the set
includes
polynomials of every non-negative degree \textit{except} those in $\Kk_N^{(\lambda)}$.

        Whenever $\lambda$ is an \emph{even} partition (which is a
        partition whose distinct entries all appear an even number
        times) and $y<0$ is any negative number then the 
        $\lambda$-XHPs are orthogonal with respect to a Hermite-like
        inner product:
\begin{align}\label{eqn:ipXHP}
    \langle H^{(\lambda)}_n, H^{(\lambda)}_m \rangle_\lambda = \int_{-\infty}^{\infty} H^{(\lambda)}_nH^{(\lambda)}_m \frac{e^{\frac{x^2}{4y}}}{\lp \tau^{(\lambda)}(x,y)\rp^2}\,dx = \nu_{n-N}(y)
\delta_{n,m}
\end{align}
where
$$
\tau^{(\lambda)}(x,y)=\frac{\text{Wr}(H_{k_\ell},\dots,H_{k_1})}{\prod_{i<j<\ell}(k_i-k_j)}
\hbox{ and }
\nu_n(y)=2(-\pi y)^{\frac{1}{2}}(-2y)^m \frac{(m+\ell)!}{\prod_{i=1}^\ell (m-m_i(\lambda))}.
$$

\begin{remark}
If the partition is not even, the $\lambda$-XHPs are still orthogonal under a modified version of
the inner product above where the only change is that the path of integration in the complex
plane is shifted in order to avoid the roots of
$\tau^{(\lambda)}(x,y)$ \cite{Veselov}.
  \end{remark}

  \subsection{The Adelic Grassmannian $Gr^{ad}$}
  In his two papers that introduced and studied the so-called ``adelic
  Grassmannian'', George Wilson provided two seemingly different but
  entirely equivalent descriptions \cite{Wilson1993,Wilson1998}.
They are each best described as a procedure for making functions of a
certain type, in one case using Wronskians and linear functionals and in the other
using matrix multiplication and rank one conditions.  Remarkably,  the set
of functions produced by each of the procedures is the same.

\subsubsection{$Gr^{ad}$ in Terms of Linear Functionals}
   For $n\in\mathbb{N}\cup\{0\}$ and $c\in \complex$ define:
\begin{align*}
    \delta_{c}^n(f(z))= \partial^n_z \lp f(z) \rp \Big|_{z=c}.
\end{align*}
to be the linear functional which differentiates $n$ times in the
variable $z$ and evaluates the function at $z=c$, which is called the
support of the functional.

Select $n$ linear combinations $D=\{d_1,\ldots,d_n\}$ of these functionals with complex
coefficients with the properties that:
\begin{itemize}
\item The set of all $n$ of them is linearly independent (i.e. no
  linear combination of them is the functional that takes the constant
  value zero).
\item Each one of them is supported at only one point in $\complex$.  
\end{itemize}
The second property assures us that there are $n$ (not necessarily
distinct numbers) $z_j\in\complex$ such that
$$
d_j=\sum_{q=0}^{m_j} \gamma_{jq} \delta^{q}_{z_j}\hbox{ for some }\gamma_{jq} \in \complex.
$$

Let the \textit{vacuum} wave function
$$
\BPsi_0(t_1,t_2,\ldots,z)=\hbox{exp}\sum_{j=1}^{\infty}
t_jz^j
$$
be a function of the infinitely many ``time variables'' $t_j$ (for
$j\in\mathbb{N}$) and $z$.

From the selection of $D$, we produce a modified wave function
according to the formula
$$
\BPsi_D(t_1,t_2,\ldots,z)=\frac{\hbox{Wr}(\phi_1,\ldots,\phi_n,\BPsi_0)}{\hbox{Wr}(\phi_1,\ldots,\phi_n)\prod_{j=1}^n(z-z_j)}
  $$
  where the Wronskian involves differentiation with respect to $t_1$
  and $\phi_j=d_j(\BPsi_0)$ is the result of applying the
  $j^{th}$ functional to the vacuum wave function.
Notice that the $z_j$'s appearing in the denominator are the supports
of the distributions in $D$.

One definition of $Gr^{ad}$ is that it is the collection of all wave
functions $\BPsi_D$ which can be made in this way using
\textit{any} selection of functionals $D$ meeting the criteria above.

\begin{remark}
  It should be noted that the map $D\mapsto \BPsi_D$ is
  \textit{not} one to one.  In particular, if $D$ and $D'$ span the
  same subspace then $\BPsi_D=\BPsi_{D'}$.
  But, the lack of uniqueness is much worse than that.  In fact, it is
  always possible to start with a choice of $D$ with $n$ elements,
  ``integrate'' the elements in a certain way, and then add
  $\delta_c^0$ to yield an alternative set containing $N+1$ elements
  corresponding to the same wave function.  However, we will not be
  elaborating further on this issue in this note.
\end{remark}

\begin{remark}\label{rmk:KPpsi}
  The wave function $\BPsi_D$ is connected to the KP
  Hierarchy in the following sense: Write
  $$
  \BPsi_D(t_1,t_2,\ldots,z)=\left(1+\sum_{j=1}^{\infty}
    \alpha_j(t_1,t_2,\ldots)z^{-j}\right)\BPsi_0.
  $$
  Let $W=1+\sum \alpha_j\partial_{t_1}^{-j}$ be the
  pseudo-differential operator in the variable $t_1$ with the same
  coefficients and let $\mathcal{L}=W\circ \partial_{t_1}\circ W^{-1}$
  be the result of conjugating the simple differential operator
  $\partial_{t_1}=\frac{d}{dt_1}$ by $W$.  Then $\mathcal{L}$
  satisfies the hierarchy of equations
  $$
  \frac{\partial}{\partial
    t_j}\mathcal{L}=[(\mathcal{L}^j)_+,\mathcal{L}]
  $$
  for every positive integer $j$ (where the positive subscript
  indicates projection onto differential operators).  Consequently,
   the coefficients of $\mathcal{L}$ satisfy integrable nonlinear
  partial differential equations, including many from mathematical
  physics such as the KP Equation that models the behavior of water
  waves. 
(For more information about these objects and their role in soliton theory, see  \cite{KasmanGOST}
or \cite{SegalWilson}.)
\end{remark}
\begin{remark}\label{rem:xt1}
 Following standard convention in the context of the KP hierarchy, this paper will identify the variables  $x$ (appearing in the orthogonal polynomials) and $t_1$ (the first time variable of the KP Hierarchy), considering them to be two different expressions representing the same variable. 

\end{remark}

\subsubsection{$Gr^{ad}$ and CM Pairs}
When the commutator of two operators is equal to the identity
operator, they are said to satisfy the ``canonical commutation
relation''.  This identity is of fundamental importance
in quantum mechanics.
However, it is easy to see that finite-dimensional matrices over $\mathbb{C}$ cannot satisfy the
canonical commutation relation\footnote{Since the properties of the trace imply that $\textup{tr}[X,Z]=0$, the commutator of two $n\times n$ matrices can only equal the identity if the scalar field has positive characteristic so that $n\cdot1=1+1+\cdots+1=0$.}.  Since the commutator of two such matrices cannot differ from the identity by the zero matrix, 
one can argue that the closest that they can get to that situation is if their commutator differs from the identity by a rank
\textit{one} matrix.  Matrices which satisfy this ``almost canonical
commutation relation'' are called \textit{Calogero-Moser Pairs} and
they arise in the study of integrable particle systems
\cite{KKSCMPAIRS}.  We define them formally by saying:

\begin{definition}
Let $X,Z \in M_{NN}\left( \complex \right)$ be constant $N\times N$ matrices.  The ordered Pair
$\left(X,Z \right)$ is said to be a Calogero-Moser Pair (CM Pair) if and only if $\text{Rank}\left( [X,Z] - I \right)=\text{Rank}\lp XZ-ZX-I \rp=1$.
\end{definition}

Let $(X,Z)$ be any CM Pair.  Then, since $[X,Z]-I$ is a rank one
matrix, there exist a row vector  $\boldsymbol{a}$ and a column vector
$\boldsymbol{b}$ such that $[X,Z]-I=\boldsymbol{b}\boldsymbol{a}$.
Wilson \cite{Wilson1998} showed that the function
$$
\BPsi^{(X,Z)}(t_1,t_2,t_3,\ldots,z)=(1+\boldsymbol{a}\tilde
  X^{-1}\tilde Z^{-1}\boldsymbol{b})\BPsi_0
  $$
  is a KP wave function in $Gr^{ad}$  (cf.~Remark~\ref{rmk:KPpsi}) where
  $$
  \tilde X=-X+\sum_{j=1}^{\infty} j t_j
  Z^{j-1}\qquad\hbox{and}\qquad\tilde Z=zI-Z.
  $$
  Moreover, running over all CM
  Pairs this produces exactly the same set of functions as one gets
  by computing $\BPsi_D$ as above and running through all
  choices of functionals $D$.  It is in that sense that $Gr^{ad}$ can
  be equivalently described in terms of functionals or CM Pairs.

  \begin{remark}
 Although    there is not a unique choice of vectors $\boldsymbol{a}$ and
    $\boldsymbol{b}$ corresponding to a fixed  $(X,Z)$, the wave
    function depends only on the choice of the CM Pair and not on
    the vectors chosen.  Moreover, the map $(X,Z)\mapsto
    \BPsi^{(X,Z)}$ is not one-to-one in the sense that for
    any invertible $N\times N$ matrix $S$, the matrices $SXS^{-1}$ and
    $SZS^{-1}$ are also a CM Pair and they correspond to the exact
    same wave function.
  \end{remark}


    \subsubsection{Stationary Wave Functions and the Bispectral Involution}
The main question of interest in Wilson's papers \cite{Wilson1993} and
\cite{Wilson1998} was the \textit{bispectrality} of the wave functions
in $Gr^{ad}$.  In that context, the infinitely many time variables
$t_i$ for $i\geq 2$ are unnecessary.  For a wave function
$\BPsi(t_1,t_2,t_3,\ldots,z)$ in $Gr^{ad}$, the result of
setting $t_1=x$ and $t_i=0$ for $i\geq 2$ is called the corresponding
stationary wave function.  Since the map
$$
\BPsi(t_1,t_2,t_3,\ldots,z)
\mapsto
\BPsi(x,0,0,\ldots,z)
$$
is one-to-one, no information is lost by considering only stationary
wave functions.  Furthermore, an interesting and important symmetry of
$Gr^{ad}$ which was discovered by Wilson is
most naturally written in terms of CM Pairs and their stationary wave
functions: 
 Note that if $(X,Z)$ is a CM Pair then so is
$(Z^{\top},X^{\top})$.  The corresponding stationary wave functions
are related by nothing more than an exchange of the two variables:
\begin{equation}\label{eqn:bispinv}
\BPsi^{(Z^{\top},X^{\top})}(x,0,0,\ldots,z)=
\BPsi^{(X,Z)}(z,0,0,\ldots,x).
\end{equation}

\subsection{XHPs Hiding in $Gr^{ad}$}
It was a surprise in 2020 when it was realized that the Exceptional
Hermite Polynomials have a generating function that is a KP wave
function from $Gr^{ad}$ \cite{MilsonKasman2020}.  In the terminology and notation of this
paper, this result can be most simply stated in the following way:

Choose any partition $\lambda$ and let $D=\{\delta_0^{k_1},\cdots,\delta_0^{k_\ell}\}$ be
the set of linear functionals supported at $z=0$ having orders
determined by the set $\Kk^{(\lambda)}=\lb k_1,\dots,k_\ell \rb$ from
Section~\ref{subsection:Partitions and Young Diagrams}. Then define
\begin{equation}
\BPsi^{(\lambda)}(x,y,z)=\BPsi_D(x,y,0,0,\ldots,z)\label{eqn:XHPgenfunc}
\end{equation}
to be the corresponding wave function with $t_1=x$, $t_2=y$ and
$t_i=0$ for $i\geq3$.   The function $\BPsi^{(\lambda)}$ is a generating
function in the sense that
\begin{align}
    {\BPsi}^{(\lambda)}(x,y,z)&=\series_{m=-\ell}\frac{\prod^\ell_{i=1}(m-m_i(\lambda))}{(m+\ell)!}\frac{H_{m+N}^{(\lambda)}(x,y)}{\tau^{(\lambda)}(x,y)}z^m
\end{align}
or equivalently
\begin{align}                                
    z^N\tau^{(\lambda)}(x,y) {\BPsi}^{(\lambda)}(x,y,z)&=\series_{n=0} \frac{\prod_{k\in \Kk_N}(n-k)}{n!}H^{(\lambda)}_n(x,y)z^n. \label{eqn:eqnforXHPgen}
\end{align}

The usefulness of this observation in the paper
\cite{MilsonKasman2020} stemmed mostly from the bispectrality of this
wave function.  Here, however, we have a different focus: what is the
CM Pair corresponding to this wave function and how can it be used to
further study the Exceptional Hermite Polynomials?

\section{The CM Pair for the $\lambda$-XHPs}
It follows from Equation~\eqref{eqn:XHPgenfunc} that the generating function $\BPsi^{(\lambda)}(x,y,z)$ of the
$\lambda$-XHPs 
 is a
  KP wave function from $Gr^{ad}$ (with the time variables set
  appropriately).  However, we also know that every KP wave function
  in $Gr^{ad}$ can be written in terms of CM Pairs.  So, the question
  arises: What CM Pair corresponds to the generating function for the
  $\lambda$-XHPs?

Before we answer that question, we introduce the following notation:
  \begin{definition}\label{def:elemmat}  The elementary matrix
  $E_{\mu\nu}(i,j)$ is defined to be the $\mu\times\nu$ matrix
  whose only non-zero entry is a $1$ in the $i^{th}$ row and $j^{th}$
  column.  If the size of the matrix is clear from context, the
  subscript may be omitted so that it is written simply as $E(i,j)$.
  \end{definition}

  Now, we are ready to give a matrix pair which answers the question,
  but it is necessary to write the partition using Frobenius
  notation (cf. Section~\ref{sec:frob}).

  \begin{definition}\label{def:CMPairForAnyPartition}
Fix a partition
$\lambda=\{(\alpha_1,\dots,\alpha_J)|(\beta_1,\dots,\beta_J)\}$ of $N$.
Let $N_k=\alpha_k+\beta_k+1$ for $1\leq k\leq J$ and define
$X_{ij}^{\lambda}$ and $Z_{ij}^{\lambda}$ to be the $N_i\times N_j$ matrices
\begin{align*}
    X_{ii}^\lambda&=-\sum_{n=1}^{\beta_i} nE(n+1,n)+\sum_{n=1}^{\alpha_i} nE(N_i+1-n,N_i-n)  \hspace{.15cm}  \\
    X_{ij}^\lambda&= N_j \sum_{n=0}^{\alpha_j} E(\beta_i+2+n,\beta_j+1+n) &\text{for}& \hspace{.15cm} &i<j \\
    X_{ij}^\lambda&= -N_j \sum_{n=0}^{\beta_i} E(\beta_i+1-n,\beta_j-n) &\text{for}& \hspace{.15cm} &j<i \\
    Z_{ii}^\lambda&=\sum_{n=1}^{N_i-1}E(n,n+1) \hspace{.15cm} \\
    Z_{ij}^\lambda&=\boldsymbol{0} &\text{for}& \hspace{.15cm} &i\neq j. \\
\end{align*}
The $N\times N$ matrices $X^\lambda$ and $Z^\lambda$ with block
decomposition
$$
X^\lambda=\left[\begin{matrix}X_{11}^\lambda&\cdots&X_{1J}^\lambda\\
    \vdots&\ddots&\vdots\\
    X_{J1}^\lambda&\cdots&X_{JJ}^\lambda\end{matrix}\right]
\qquad
\hbox{and}
\qquad
Z^\lambda=\left[\begin{matrix}Z_{11}^\lambda&\cdots&Z_{1J}^\lambda\\
    \vdots&\ddots&\vdots\\
    Z_{J1}^\lambda&\cdots&Z_{JJ}^\lambda\end{matrix}\right]
$$
are called the CM Pair associated to $\lambda$.
\end{definition}

Of course, it is not yet clear that these matrices
are related in any useful way to the $\lambda$-XHPs, or even that they
are actually a CM Pair.  The next two results from \cite{Thesis} show
that they are:

\begin{theorem}\label{thm:theyareCMPair}
The matrices $X^\lambda$ and $Z^\lambda$ satisfy the rank one condition
  $$
  [X^\lambda,Z^\lambda]-I=\boldsymbol{b}\boldsymbol{a}
  $$
where $\boldsymbol{b}$ and $\boldsymbol{a}$ are the $N \times 1$ and
$1 \times N$ matrices with block decompositions
\begin{align}
    \boldsymbol{b}&=\begin{bmatrix}
    E_{N_1,1}(\beta_1+1,1) \\
    \vdots \\
    E_{N_i,1}(\beta_i+1,1) \\
    \vdots \\
    E_{N_J,1}(\beta_J+1,1)
    \end{bmatrix}
  \hspace{.1cm}&\boldsymbol{a}^\top= \begin{bmatrix}
    -N_1E_{N_1,1}(\beta_1+1,1) \\
    \vdots \\
    -N_jE_{N_j,1}(\beta_j+1,1) \\
    \vdots \\
    -N_JE_{N_J,1}(\beta_J+1,1)
    \end{bmatrix}
\label{eqn:intuitiveFormofBA}.
\end{align}
  and hence   $(X^\lambda,Z^\lambda)$
is a CM Pair.
\end{theorem}

Furthermore, not only is $(X^{\lambda},Z^{\lambda})$ \textit{some} CM
Pair, it is one in terms of which the generating
function can be written:
\begin{theorem}\label{thm:rightone}
The generating function $\BPsi^{(\lambda)}(x,y,z)$ of the
$\lambda$-XHPs can be written in terms of the CM Pair
$(X^\lambda,Z^\lambda)$ as
\begin{equation}
\BPsi^{(\lambda)}(x,y,z)=\BPsi^{(X^\lambda,Z^\lambda)}(x,y,0,0,\ldots,z)
=(1+\boldsymbol{a}\tilde
  X^{-1}\tilde Z^{-1}\boldsymbol{b})\BPsi_0
\label{eqn:XHPgenfuncCM}
\end{equation}
with $\boldsymbol{a}$ and $\boldsymbol{b}$ as in
\eqref{eqn:intuitiveFormofBA}
and where
  $$
  \tilde X=-X^\lambda+xI+2yZ^{\lambda}\qquad\hbox{and}\qquad\tilde Z=zI-Z^\lambda.
  $$
\end{theorem}
\begin{remark}
  The proof of Theorem~\ref{thm:theyareCMPair} 
is
  nothing more than a long and tedious computation.  The proof of
  Theorem~\ref{thm:rightone} given in \cite{Thesis}, on the other hand, involves some
  interesting ideas.  In particular, it utilizes the fact that wave
  functions in $Gr^{ad}$ can be identified by the finitely-supported
  linear functionals which annihilate them.  It is relatively easy to
  see that in the case of $\BPsi_D$ those would be precisely the ones
  which multiply by $\prod(z-z_j)$ and then apply some linear
  combination of the functionals in $D$.  On the other hand, Lemma~3.1
  from the paper \cite{BGK} provides a method involving residues for
  finding the linear functionals which annihilate a wave function
  $\BPsi^{(X,Z)}$ written in terms of a CM Pair.  Confirming that
  these collections of functionals are equal in the special case $D=D_\lambda$ and
  $(X,Z)=(X^\lambda,Z^\lambda)$ was the proof method utilized in \cite{Thesis}.  (Another
   proof method of an equivalent result can be found in the appendix by Macdonald to
  Wilson's 1998 paper \cite{Wilson1998}.  We became aware of this
  alternate proof only afterwards.)
\end{remark}

\begin{example}\label{exampI}
  We will illustrate this procedure for associating a CM Pair to a partition with an example. Choose the partition $\lambda=\lb 2,2,1,1 \rb$. Using $(\ref{eqn:MayaDiagramDef})$, $(\ref{eqn:indexSet})$, and $(\ref{eqn:scriptI})$ we can conclude that:
$$
        \Kk_N^{(\lambda)}=\lb 7,6,4,3,1,0 \rb\qquad\hbox{ and } \qquad\Ii=\lb 2,5,8,9,\dots \rb.
$$
    Now we use Definition \ref{def:CMPairForAnyPartition} in order to construct:
\begin{align*}
    &X=\begin{bmatrix}
       0 & 0 & 0 &0 &0 &0 \\
       -1 & 0 & 0 &0 &0 &0 \\
       0 & -2 & 0 &0 &0 &0 \\
       0 & 0 & -3 &0 &0 &0 \\
       0 & 0 & 0 &1 &0 &1 \\
       0 & 0 & -5 &0 &0 &0
    \end{bmatrix}& &Z=\begin{bmatrix}
        0 & 1 & 0 &0 &0 &0 \\ 
        0 & 0 & 1 &0 &0 &0 \\ 
        0 & 0 & 0 &1 &0 &0 \\ 
        0 & 0 & 0 &0 &1 &0 \\ 
        0 & 0 & 0 &0 &0 &0 \\ 
        0 & 0 & 0 &0 &0 &0 \\ 
    \end{bmatrix}.&
\end{align*}
Then $[X,Z]-I=\boldsymbol{b}\boldsymbol{a}$ where:
$$
    \boldsymbol{b}^\top=\begin{bmatrix}
        0 & 0 &0 &1 & 0 &1 
    \end{bmatrix}\hbox{ and }\boldsymbol{a}=\begin{bmatrix}
        0 & 0 &0 &-5 & 0 &-1 
    \end{bmatrix}.
$$
Now we can form 
\begin{align*}
&\tilde{X}=-X+xI+2yZ=\begin{bmatrix}
       x & 2y & 0 &0 &0 &0 \\
       1 & x & 2y &0 &0 &0 \\
       0 & 2 & x &2y &0 &0 \\
       0 & 0 & 3 &x &2y &0 \\
       0 & 0 & 0 &-1 &x &-1 \\
       0 & 0 & 5 &0 &0 &x
    \end{bmatrix}.& 
    && &&
\end{align*}
The corresponding KP wave function is then given by $(\ref{eqn:XHPgenfuncCM})$.
This computation will be continued below in Example~\ref{exampII} where the  exceptional orthogonal polynomials associated to the partition $\lambda$ will be produced directly from these matrices.
\end{example}

\section{Formula for $\lambda$-XHPs in terms of their CM Pair}
The main result of the thesis is the following formulas which give any
Exceptional Hermite Polynomial as a linear combination of classical
Hermite Polynomials with coefficients written in terms of the
corresponding CM Pair:
\begin{theorem} \label{thm:XHPsExpanded}
  Let $\lambda$ be a partition of $N$, $(X,Z)=
  (X^{\lambda},Z^{\lambda})$ from
Definition~\ref{def:CMPairForAnyPartition}, $\tilde
X=-X+xI+2yZ$, and
$\boldsymbol{a},\boldsymbol{b}$ be the vectors from \eqref{eqn:intuitiveFormofBA}.
The $\lambda$-XHPs $\{H_n^{(\lambda)}(x,y)\ :\
n\in\Ii^{(\lambda)}\}$
can be computed using the formulas:
{  \begin{align*}
    &H^{(\lambda)}_n= \frac{\det(\tilde X)}{s(n)}\left( \sum_{k=0}^n \frac{n!}{ (n-k)! } H_{n-k}(x,y) \boldsymbol{a}\tilde{X}^{-1}Z^{N-k-1}\boldsymbol{b} \right)& \\
    &H^{(\lambda)}_n= \frac{\det(\tilde X)}{s(n)}\left(
      \sum_{k=0}^{N-1} \frac{n!}{ (n-k)!}H_{n-k}(x,y)\boldsymbol{a}
      \tilde{X}^{-1}Z^{N-k-1}\boldsymbol{b}+ \frac{n!}{(n-N)!}
      H_{n-N}(x,y) \right)& 
\end{align*}}
when $n<N$ and $n \geq N$ respectively and where $s(n)=\prod_{k\in \Kk^{(\lambda)}_N}(n-k)$.
\end{theorem}

\begin{remark}
  The proof of Theorem~\ref{thm:XHPsExpanded} begins with the
  observation that applying $\delta^n_0$ to $(\ref{eqn:eqnforXHPgen})$ yields
\begin{align}
    \delta_{0}^n \lp z^N \tau^{(\lambda)}(x,y) \BPsi^{(\lambda)} \rp  = s(n)H_n^{(\lambda)}(x,y) \label{eqn:almostDefforXHPS}.
\end{align}
  Since Equation
\eqref{eqn:XHPgenfuncCM} gives the generating function in terms of
matrix products, it is relatively easy to apply the functional
$\delta_0^n$ to it.  In addition, two identities were used to further
simplify the result:  The rank one condition satisfied by the CM Pair
implies that
$$
ZX+\boldsymbol{b}\boldsymbol{a}=X Z-I
$$
and the Matrix Determinant Lemma \cite{wikiMDL} further allows one to
translate back and forth between matrix products and determinants.
Additionally, we made use of the correspondence that
$\tau^{(\lambda)}(x,y)=\det(\tilde X)$.
\end{remark}
\begin{example}\label{exampII}
Recall Example~\ref{exampI}
where specific matrices $X$, $Z$, $\tilde{X}$, $\boldsymbol{a}$, and $\boldsymbol{b}$
were associated to the choice of partition $\lambda=\lb 2,2,1,1 \rb$.
Here we will illustrate the use of those matrices in producing the corresponding orthogonal polynomials 
$H_n^{(\lambda)}(x,y)$. According to Theorem
\ref{thm:XHPsExpanded} we  can calculate $H_n^{(\lambda)}(x,y)$ using matrix products and the determinant of $\Tilde{X}$ for
all $n\in \Ii$.  Here, for example, are the cases $n=2$, $5$ and $8$:
\begin{align*}
	H^{(\lambda)}_2(x,y) &= \frac{\text{det}\lp \Tilde{X} \rp}{80}\lp \sum_{k=0}^2 \frac{2!}{(2-k)!} H_{2-k}(x,y) \boldsymbol{a}\Tilde{X}^{-1}Z^{6-k-1}\boldsymbol{b} \rp \\
	&=\frac{\text{det}\lp \Tilde{X} \rp}{80}\lp 2 H_0(x,y) \boldsymbol{a}\Tilde{X}^{-1}Z^{3}\boldsymbol{b}  \rp \\
	&=x^2-2y \\
	H^{(\lambda)}_5(x,y) &= \frac{\text{det}\lp \Tilde{X} \rp}{80}\lp \sum_{k=0}^5 \frac{5!}{(5-k)!} H_{5-k}(x,y) \boldsymbol{a}\Tilde{X}^{-1}Z^{6-k-1}\boldsymbol{b} \rp \\
	&=\frac{\text{det}\lp \Tilde{X} \rp}{80} ( 20H_3(x,y) \boldsymbol{a}\Tilde{X}^{-1}Z^{3}\boldsymbol{b}+60H_2(x,y) \boldsymbol{a}\Tilde{X}^{-1}Z^{2}\boldsymbol{b} \\
	&\;\;\;+ 120H_1(x,y) \boldsymbol{a}\Tilde{X}^{-1}Z\boldsymbol{b}+120H_0(x,y)\boldsymbol{a}\Tilde{X}^{-1}I\boldsymbol{b} ) \\
                             &=x^5-20x^3y+60xy^2\\
      H^{(\lambda)}_8(x,y)&=\frac{\det(\tilde
                            X)}{2240}\left(\sum_{k=0}^5\frac{8!}{(8-k)!}H_{8-k}(x,y)\boldsymbol{a}\tilde
                            X^{-1}Z^{5-k}\boldsymbol{b}+\frac{8!}{2!}H_2(x,y)\right)\\
  &= x^8-20x^6y+120x^4y^2-240x^2y^3-240y^4
\end{align*}
\end{example}

\section{Finitely-Supported Distributions Annihilating $\lambda$-XHPs}

For $n\in\nats$ and $c\in\complex$, let $\Delta_c^n$ be the linear
functional 
$$
\Delta_c^n(f(x))=f^{(n)}(c)
$$
which differentiates $n$ times in the variable $x$ and then evaluates
at $x=c$.  From an abstract perspective, this is the same functional
as $\delta_c^n$ introduced earlier, with the only difference being
that the this new one acts on functions of the variable $x$ while the
one we saw earlier acted on functions of $z$.  (Since we will be applying these finitely supported distributions to functions that depend on both variables, this notation allows us to indicate on which variable they are acting.  For example, if $\psi(x,z)=x^2z^3$ then $\Delta_1^1(\psi(x,z))=2z^3$ while $\delta_1^1(\psi(x,z))=3x^2$.)

For any fixed partition $\lambda$ it is known that there are
functionals which can be written as linear combinations of the
$\Delta_c^n$ that annihilate the entire family of $\lambda$-XHPs.
In particular, \cite{MilsonKasman2020} explains how one can find them
using Wilson's bispectral involution on $Gr^{ad}$.  However, the bispectral
involution is not defined in an explicit form when $Gr^{ad}$ is
described using functionals in $z$ and so \cite{MilsonKasman2020} was
only able to give an implicit description.

In contrast, the bispectral involution is not only explicitly
described but also very simple in terms of CM Pairs.  Combining
\eqref{eqn:bispinv} together with Theorem~\ref{thm:rightone} allowed
us to prove the following result giving formulas for those functionals
annihilating the $\lambda$-XHPs for almost all values of $y$:

\begin{theorem}\label{thm:annh}
Let $\lambda$ be a partition of $N$ and let $(X,Z)=(X^\lambda,Z^\lambda)$ be
the corresponding CM Pair as in
Definition~\ref{def:CMPairForAnyPartition}.  Now, suppose $y$ is
chosen so that the matrix $X-2yZ$ has distinct
  eigenvalues $\{\gamma_1,\ldots,\gamma_N\}$.   Then, let $S$ be a
  matrix such that $S^{-1}(X^\top-2yZ^\top)S=\textup{diag}(\gamma_1,\ldots,\gamma_N)$ and define
  $$
  X'=S^{-1}Z^\top S,\ Z'=\textup{diag}(\gamma_1,\ldots,\gamma_N),\
  \boldsymbol{b}'=S^{-1}\boldsymbol{a}^{\top},\ \hbox{and}\
  \boldsymbol{a}'=\boldsymbol{b}^\top S
$$
where $\boldsymbol{a}$ and $\boldsymbol{b}$ are as in \eqref{eqn:intuitiveFormofBA}.
Finally, for $1\leq i\leq N$ define the functions
\begin{equation*}
	g_i(x)=\frac{\boldsymbol{a}'_{i}\boldsymbol{b}'_{i}+X'_{ii}(x-\gamma_i)}{\prod_{n\neq i}\lp x-\gamma_n \rp}.\label{eqn:gsubi}
\end{equation*}
Then for $1\leq i\leq N$ the functionals
$$\theta_i=g_i(\gamma_i)\Delta^1_{\gamma_i}+g_i'(\gamma_i)\Delta^0_{\gamma_i}$$
annihilate each of the $\lambda$-XHPs.  Furthermore, a
finitely-supported functional $d$ satisfies $d(H_n^{(\lambda)})=0$ for
all $n\in\Ii^{(\lambda)}$ if and only if $d\in\textup{span}(\theta_1,\ldots,\theta_N)$.
\end{theorem}

\begin{example}
Consider the partition $\lambda=\{1,1\}=\{(0)|(1)\}$.  The
corresponding CM pair is
$$
X=\left[
\begin{array}{cc}
 0 & 0 \\
 -1 & 0 \\
\end{array}
\right]\qquad Z=\left[
\begin{array}{cc}
 0 & 1 \\
 0 & 0 \\
\end{array}
\right]
$$
which satisfy
$$
[X,Z]-I=\left[
\begin{array}{cc}
 0 & 0 \\
 0 & -2 \\
\end{array}
\right]=\boldsymbol{b}\boldsymbol{a}=\left[
\begin{array}{c}
 0 \\
 1 \\
\end{array}
\right]\left[
\begin{array}{cc}
 0 & -2 \\
\end{array}
\right].
$$
From these one can make the KP wave function
$$\psi^{(X,Z)}(x,y,z)=e^{x z+y z^2} \left(\frac{2}{z^2 \left(x^2-2
      y\right)}-\frac{2 x}{z \left(x^2-2 y\right)}+1\right)
$$
with corresponding $\tau$-function
$$
\tau(x,y)=x^2-2y.
$$
Clearing the denominator of the wave function and expanding as a
series in $z$ yields:
\begin{eqnarray*}
\tau(x,y)z^2\psi(x,y,z)&=&2+0z+0z^2+\frac{1}{3} \left(x^3-6 x y\right)  z^3+\frac{1}{4}
   \left(x^4-4 x^2 y-4 y^2\right)  z^4\\&&+\frac{1}{10} 
   \left(x^5-20 x y^2\right) z^5+O\left(z^6\right).
   \end{eqnarray*}
  The non-zero coefficients in that expression are
 (up to constant scalar multiples)   the family of Exceptional Hermite Polynomials corresponding to
   $\lambda$.

   This section explains how we can use the CM pair to produce differential conditions in
   the variable $x$ which are satisfied by all of the polynomials in
   that sequence whenever $y$ is chosen so that $\tau(x,y)$ has
   distinct roots.  Allow us to demonstrate this procedure.

   To keep things simple, let us begin by choosing $y=-2$ so that
   $$
   \tau(x,-2)=x^2+4=(x-2\complexi)(x+2\complexi)
   $$
   has the distinct roots $\pm 2\complexi\in\complex$.  It is no
   coincidence that these same numbers are the eigenvalues of the
   matrix $X-2yZ$ (and its transpose)  We will focus here only on this
   value of $y$ and on the eigenvalue $2\complexi$, though the
   procedure works the same for any $y$ for which the eigenvalues are
   distinct and for any one of those eigenvalues.

  It is easy to see that
   $$
   X^\top-2yZ^\top=\left[
\begin{array}{cc}
 0 & -1 \\
 -4 \complexi & 0 \\
\end{array}
\right]
$$
and we know that this matrix is diagonalizable (because of how $y$ was
chosen).  In fact,
$$
Z'=S^{-1}(X^\top-2yZ^\top)S=\hbox{diag}(-2\complexi,2\complexi)
\hbox{ where }
S=\left[
\begin{array}{cc}
 -\complexi & \complexi \\
 2 & 2 \\
\end{array}
\right].
$$
We use that same matrix $S$ to define
$$
X'=S^{-1}Z^\top S,\ \boldsymbol{a}'=\boldsymbol{b}^\top S, \hbox{ and
}\boldsymbol{b}'=S^{-1}\boldsymbol{a}^\top.
$$
They turn out to be:
$$
X'=\left[
\begin{array}{cc}
 -\frac{\complexi}{4} & \frac{\complexi}{4} \\
 -\frac{\complexi}{4} & \frac{\complexi}{4} \\
\end{array}
\right]
,\
\boldsymbol{a}'=\left[
\begin{array}{cc}
 2 & 2 \\
\end{array}
\right]
,\ \hbox{and}\
\boldsymbol{b}'=
\left[
\begin{array}{c}
 -\frac{1}{2} \\
 -\frac{1}{2} \\
\end{array}
\right]
.
$$
These form a new CM pair because
$$
[X',Z']-I=\boldsymbol{b}'\boldsymbol{a}'.
$$
(Essentially, what we have done is shift the pair $(X,Z)$ using the
second flow of the KP hierarchy and then applied the bispectral
involution).

These matrices give us the values we need to construct those
differential conditions satisfied by the $\lambda$-XHPs.  Since the
eigenvalue of interest here is $2\complexi$ which is in the second row
and column of $Z'$, we will use $i=2$ in the formulas.  (If instead we
were looking at the eigenvalue $-2\complexi$ we would have been using
$i=1$.)  In particular,
we compute the function $g_i(x)$ from Equation~\eqref{eqn:gsubi}
\begin{eqnarray*}
g_2(x)&=&\frac{\boldsymbol{a}_2'\boldsymbol{b}_2'+X'_{22}(x-2\complexi)}{x+2\complexi}
=\frac{-1+\frac{1}{4} \complexi (x-2 \complexi)}{x+2 \complexi}
\end{eqnarray*}
(where the subscripts on $\boldsymbol{a}'$, $\boldsymbol{b}'$ and $X'$
just indicate that we should select the corresponding entry from the
constant matrices of those names listed above).

The theorem then tells us that the linear functional
$$
g_2(2\complexi)\Delta_{2\complexi}^1+g_2'(2\complexi)\Delta_{2\complexi}^0=\frac{\complexi}{4}\Delta_{2\complexi}^1$$
will annihilate all of the $\lambda$-XHPs.  Indeed, one can check that
each of the coefficient functions in the series expansion
above has the property that its derivative with
respect to $x$ evaluated at $x=2\complexi$ is zero.
\end{example}

\section{Concluding Remarks}
\subsection{Relationship to the Thesis}
This opportunity to report on the main results in \cite{Thesis} also
has given us a chance to correct some of the errors that thesis
contained.  
For example, the definition of $g_i(x)$ given in
Theorem~\ref{thm:annh}  was written incorrectly in the thesis, with $\gamma_i$ being written in place of $X'_{ii}$ by mistake.
 Furthermore, the
statement of Theorem~\ref{thm:XHPsExpanded} in the thesis
unfortunately confused $n$ and $N$.  In general, wherever statements
in this
research announcement differ from those in the thesis, it should be
understood that a small error has been corrected.

On the other hand, the proof methods utilized in \cite{Thesis} remain
essentially valid (with those small corrections implemented), and so
one should still consult that original document for more details on
the proofs.

The thesis also falsely claims that it is \textit{always} possible to
choose $y$ so that the roots of $\tau^{(\lambda)}(x,y)$ (equivalently
the eigenvalues of $X-2yZ$) are distinct.  If the corresponding
solution to the KP Hierarchy happens to be a \textit{KdV} solution
then $\tau^{(\lambda)}(x,y)$ will be independent of $y$ and the only
roots will be $x=0$.

\subsection{Future Outlook}
This is the first investigation about the connections between
Exceptional Hermite Polynomials and Calogero-Moser Pairs.
The results above conclusively demonstrate that this correspondence
can be fruitful.  However, there is more to do.

One of the main goals we were unable to achieve was to
create an alternate proof for orthogonality of the XHPs using CM Pairs. As stated in
Section \ref{subsection:XHPS}, it is already well known that the
$\lambda$-XHPs are orthogonal with respect to the inner product \eqref{eqn:ipXHP}. However, the proof of this fact requires some higher level analytic techniques. As a result we felt it would be fruitful to give a proof using Theorem \ref{thm:XHPsExpanded}, elementary linear algebra, and some basis properties of integrals.
Unfortunately we were not able to devise such a proof. 

There are many CM Pairs other than the special ones
$(X^\lambda,Z^\lambda)$ corresponding to partitions that appear in all
of these results.  There may be a way to generalize these results to a
larger subset of them, or even perhaps to all of $Gr^{ad}$, though we
did not find anything useful to say in that more general case.

It would be nice to know whether or not Theorem \ref{thm:annh} works for all $y \neq 0$. As is currently
stated, we only have that it works everywhere except on a set of
Lebesgue measure $0$. The problem amounts to showing whether or not
the matrix $\displaystyle{-X+2yZ}$ is diagonalizable.  Furthermore,
it ought to be possible to generalize Theorem~\ref{thm:annh} so that
it works also in the case that the eigenvalues are \textit{not}
distinct.  The same method should work, but it gets very complicated.

\par\medskip\par

\noindent\textit{Acknowledgements:} Both authors thank the organizers
for the opportunity to present at this fantastic special session and
the participants for their talks and feedback.  The first
author would also like to thank the other members of his thesis committee
-- Annalisa Calini and Rob Milson -- as well as Dan Maroncelli for their
assistance and support.  Finally, we are grateful to the anonymous referees for their helpful suggestions.

\bibliographystyle{amsplain}
\bibliography{kasman-conm.bib}

\end{document}